\documentclass{article}
\usepackage{amssymb}
\usepackage{amsmath}
\title{Geodesic Paths in the Finite Dimensional Unit Sphere under Sup Norm  }
\author{Teck-Cheong Lim \\ Department of Mathematical Sciences \\ George Mason University \\ 4400, University Drive\\Fairfax, VA 22030\\U.S.A.\\
{\it e-mail address}: tlim@gmu.edu}
\date{ }

\newtheorem {theo} {Theorem}

\newtheorem {rem} {Remark}
\newtheorem {co} {Corollary}

\begin{document}
\maketitle

\begin{abstract}
By geodesic path between two points $A, B$, we mean a curve whose arc length is the shortest among all curves on the surface joining $A$ and $B$. Generally there are more than one such path. Geodesic paths in the unit sphere $\{x: \|x\|_{\infty}=1\}$ of $\mathbb{R}^n$ with the sup norm are investigated.
\end{abstract}

{\it Keywords}: Geodesic, sup norm, unit sphere (cube) \\ \\
Let $X$ be the unit sphere of $\mathbb{R}^n$ with the sup norm $\|\cdot\|_\infty$ ( the $n$-dimensional $\ell^\infty$), i.e. \[X=\{x\in \mathbb{R}^n: \|x\|_\infty =1\}.\]
Let the distance $d$ on $X$ be the geodesic distance induced by the sup norm.\\
Given any two points $A, B$ on $(X,d)$, a geodesic path from $A$ to $B$ is a distance preserving map $\phi$ from a closed interval $[a,b]$ to $(X,d)$  such that $\phi(a)=A, \phi(b)=B$. We will make no distinction between the map and its image $\phi([a,b])$. Discussion of geodesic paths in other settings may be found in [1].\\
We will consider only paths whose nonvoid intersection with each $(n-1)$-dimensional face is a line segment, so every path from $A$ to $B$ can be written as a finite sequence of line segments $AC_1, C_1C_2,\cdots, C_nB$, where each $C_i$ lies in an $(n-2)$-dimensional face. A path $AC_1,C_1C_2,\cdots, C_nB$ is called {\em planar} if it is contained in a $n-1$-dimensional affine subspace. In the  absence of  ambiguity, if $B,C$ are on the same face, we also write $BC$ to stand for the length of the line segment joining $B,C$.

\begin{theo}\rm Let $n=3$. Let $A=(1,a_y,a_z)\in X, B=(b_x,1,b_z)\in X$. Then the geodesic distance between $A, B$ is the minimum of the following three quantities: \[\alpha:=\max\{ 2-a_y-b_x, |a_z-b_z|\},\] \[\beta:=\max\{2-a_z-b_x, 2-a_y-b_z\},\] \[\gamma:=\max\{2+a_z-b_x, 2-a_y+b_z\}.\]
The minimum is $\alpha$ if and only if at least one of the following four conditions is satisfied:
\begin{eqnarray}& &|a_z|\le a_y\\
& &|b_z|\le b_x\\
& &|a_y|\le a_z, \mbox{ and } |b_x|\le -b_z\\
& &|a_y|\le -a_z, \mbox{ and } |b_x|\le b_z
\end{eqnarray}
The minimum is $\beta$ if and only if at least one of the following four conditions is satisfied:
\begin{equation} a_y\le a_z, \mbox{ and }b_x\le b_z, \mbox{ and }a_z\ge 0, \mbox{ and } b_x\le a_y+a_z+b_z\end{equation}
\begin{equation}a_y\le a_z, \mbox{ and }b_x\le b_z, \mbox{ and }b_z\ge 0, \mbox{ and } a_y\le b_x+a_z+b_z\end{equation}
\begin{equation} a_y=a_z=1\end{equation}
\begin{equation} b_x=b_z=1\end{equation}
The minimum is $\gamma$ if and only if at least one of the following four conditions is satisfied:
\begin{equation} a_y\le -a_z,\mbox{ and } b_x\le -b_z, \mbox{ and }a_z\le 0, \mbox{ and } b_x\le a_y-a_z-b_z\end{equation}
\begin{equation}a_y\le -a_z, \mbox{ and }b_x\le -b_z, \mbox{ and }b_z\le 0, \mbox{ and } a_y\le b_x-a_z-b_z\end{equation}
\begin{equation} a_y=-a_z=1\end{equation}
\begin{equation} b_x=-b_z=1\end{equation}
Moreover, there exists a planar minimal path joining $A$ and $B$.
\end{theo}
{\bf Proof.}\\
First we prove that the minimum of \[f(z):=AC+CB=\max\{1-a_y,|z-a_z|\}+\max\{1-b_x,|z-b_z|\},\] where $C$ is of the form $(1,1,z),|z|\le 1$ is $\alpha:=\max\{ 2-a_y-b_x, |a_z-b_z|\}$. Observe that $AC+CB\ge \alpha$ for all $|z|\le 1$. We want to prove that the value $\alpha$ is attained at some $z$. \\ \\
Case (i): $|a_z-b_z|\ge 2-a_y-b_x.$ Then either $a_z-b_z\ge 2-a_y-b_x$ or \\$b_z-a_z\ge 2-a_y-b_x$. In the former case, $f$ attains the value $\alpha$ at any $z$ in the interval $[1-b_x+b_z, a_y+a_z-1]$, while in the latter case, at any $z$ in the interval $[a_z-a_y+1,b_x+b_z-1]$.\\ \\
Case (ii): $|a_z-b_z|\le 2-a_y-b_x$. Then $a_y+a_z-1\le -b_x+b_z+1$ and $b_x+b_z-1\le -a_y+a_z+1$. Since $a_y+a_z-1\le -a_y+a_z+1$ we see that  \[t_1:=\max\{a_y+a_z-1, b_x+b_z-1,-1\}\le t_2:=\min \{-a_y+a_z+1, -b_x+b_z+1,1\}.\] Then $-1\le t_1\le 1,-1\le t_2\le 1$. Direct checking shows that for any $z$ in the interval $[t_1, t_2]$, $|z-a_z|\le 1-a_y$ and $|z-b_z|\le 1-b_x$ and $AC+CB=2-a_y-b_x$, as required.\\ \\
Next we consider the path $AC_1, C_1C_2, C_2B$ where $C1:=(1,y,1), C2:=(x,1,1)$. The minimum of the total length function  \[g(x,y):=\max\{|y-a_y|,1-a_z\}+\max\{1-x,1-y\}+\max\{|x-b_x|,1-b_z\}\] over $(x,y)\in [-1,1]\times [-1,1]$ is $\beta_1:=\max\{2-a_z-b_x, 2-a_y-b_z, 2-a_z-b_z\}$. To prove this, first we observe from the definition and the  triangle inequality that $g(x,y)\ge \beta_1$ for all $x,y$. We shall show that $g$ attains this value $\beta_1$.\\
Write $x_1:=2-a_z-b_x, y_1:=2-a_y-b_z, z_1:=2-a_z-b_z$. \\
Case (a): $x_1\ge \max(y_1,z_1)$, i.e.$\beta_1=x_1$. Then we have $a_y+b_z\ge a_z+b_x$ and $b_z\ge b_x$. Note that $g(x,y)=x_1$ if and only if $|y-a_y|\le 1-a_z, 1-x\ge 1-y$ and $x-b_x\ge 1-b_z$. $(x,y)$ in $[-1,1]\times [-1,1]$ satisfy these conditions if and only if
\begin{eqnarray*} & & 1+b_x-b_z\le x\le  \min (1,1+a_y-a_z)\\
& & \max (a_y+a_z-1,x)\le y \le \min (1, 1+a_y-a_z)
\end{eqnarray*}
(Conditions $a_y+b_z\ge a_z+b_x$ and $b_z\ge b_x$ imply that such $(x,y)$ exist.)\\
Case (b): $y_1\ge \max(x_1,z_1)$, i.e.$\beta_1=y_1$. Then we have $a_z+b_x\ge a_y+b_z$ and $a_z\ge a_y$. Note that $g(x,y)=y_1$ if and only if $|x-b_x|\le 1-b_z, 1-y\ge 1-x$ and $y-a_y\ge 1-a_z$. $(x,y)$ in $[-1,1]\times [-1,1]$ satisfy these conditions if and only if
\begin{eqnarray*} & & 1+a_y-a_z\le y\le  \min (1,1+b_x-b_z)\\
& & \max (b_x+b_z-1,y)\le x \le \min (1, 1+b_x-b_z)
\end{eqnarray*}
(Conditions $a_z+b_x\ge a_y+b_z$ and $a_z\ge a_y$ imply that such $(x,y)$ exist.)\\
Case (c): $z_1\ge \max(x_1,y_1)$, i.e.$\beta_1=z_1$. Then we have $a_y\ge a_z$ and $b_x\ge b_z$. Note that $g(x,y)=z_1$ if and only if $x=y=1$.\\ \\
 For paths $AD_1, D_1D_2, D_2B$ where $D1:=(1,y,-1), D2:=(x,1,-1)$. The minimum of the total length function  \[h(x,y):=\max\{|y-a_y|,1+a_z\}+\max\{1-x,1-y\}+\max\{|x-b_x|,1+b_z\}\] over $(x,y)\in [-1,1]\times [-1,1]$ is $\gamma_1:=\max\{2+a_z-b_x, 2-a_y+b_z, 2+a_z+b_z\}$. To prove, we just replace $a_z$ by $-a_z$, and $b_z$ by $-b_z$ throughout the preceding proof.\\

In what follows we will write $a=a_y, b=a_z, c=b_x, d=b_z$. Since $|b-d|\le 2-b-d, 2+b+d$ for any number $b,d\in [-1,1]$, we see that  $\alpha\le \beta_1$ if and only if
$2-a-c\le 2-b-c$ or $2-a-c\le 2-a-d$ or $2-a-c\le 2-b-d$, which simplifies to $b\le a$ or  $d\le c$ or $b+d\le a+c$. And similarly $\alpha\le\gamma_1$ if and only if $-a\le b$ or $-c\le d$ or $-a-c\le b+d$. Therefore $\alpha=\min\{\alpha, \beta_1, \gamma_1\}$ if and only if at least one of the nine conditions below is satisfied:
\begin{eqnarray}& &b\le a \mbox{ and }-a\le b\\
& &b\le a\mbox{ and }-c\le d\\
& & b\le a\mbox{ and }-a-c\le b+d\\
& & d\le c  \mbox{ and }-a\le b\\
& & d\le c  \mbox{ and }-c\le d\\
& & d\le c \mbox{ and }-a-c\le b+d\\
& & b+d\le a+c \mbox{ and }-a\le b\\
& & b+d\le a+c \mbox{ and }-c\le d\\
& & b+d\le a+c\mbox{ and }-a-c\le b+d
\end{eqnarray}
Now we show that conditions (13)-(21) are equivalent to conditions (1)-(4). Clearly (1)$\Leftrightarrow$ (13); (2)$\Leftrightarrow$ (17); (3)$\Rightarrow$ (16); and (4)$\Rightarrow$ (14). On the other hand, (13)$\Leftrightarrow$ (1); (14)$\Rightarrow$ (1) or (2) or (4); (16)$\Rightarrow$ (1) or (2) or (3); (17)$\Leftrightarrow$ (2). Next we show that  (21)$\Rightarrow$ (1),(2),(3) or (4). To see this, consider four cases: $|b|\le a, |b|\le -a, |a|\le b, |a|\le -b$. In the first case, condition (1) is satisfied. In the second case, we have $|d|-|b|\le |b+d|\le a+c,$ so $|d|\le |b|+a+c\le c$ and condition (2) is satisfied. In the third case,  $b+d\le |b+d|\le a+c\le b+c \Rightarrow d\le c$, so condition (2) is satisfied if $d\ge -c$; if $d\le -c$, then $|c|\le -d$ and condition (3) is satisfied . Similarly in the fourth case, using $-b-d\le |b+d|\le a+c\le -b+c$ we see that  condition (2) is satisfied if $d\le c$, and condition (4) is satisfied if $d\ge c$. Next we show that (15)$\Rightarrow$ (1) or (21) or (4) and hence $\Rightarrow$ (1) or (2) or (3) or (4) by what was proved above for (21). Since generally $x\le y \Leftrightarrow |x|\le y \mbox{ or }|y|\le -x$, clearly (15)$\Rightarrow$ (1) or (21) or the condition $|a|\le -b$ and $|a+c|\le b+d$. But the latter condition implies $a-d\le -b-d\le a+c\le b+d\le a+d$, so $|c|\le d$ and condition (4) follows. Similarly (18)$\Rightarrow$ (2) or (21) or  (3); (19)$\Rightarrow$ (1) or (21) or  (3);(20) implies (2) or (21) or (4).\\ \\
Next we prove that $\beta_1\le \alpha$ if and only if either (i) $a\le b$ and $c\le d$ or (ii) $a=b=1$ or (iii) $c=d=1$. Since $\max\{b,d\}=\frac{1}{2}(|b-d|+b+d)$, we see that $|b-d|\le 2-b-d$ for all $b,d\in [-1,1]$ and equality holds if and only if $b=1$ or $d=1$. If $a\le b$ and $c\le d$, then $2-a-d\le 2-a-c$, $2-b-c\le 2-a-c$ and $2-b-d\le 2-a-c$ and hence $\beta_1\le \alpha$. If $a=b=1$ or $c=d=1$, then clearly $\beta_1=\alpha$. This proves the sufficiency of (i)-(iii). Next assume that $\beta_1\le \alpha$. Suppose $b<a$. If $\alpha=2-a-c$, then $2-b-c\le \beta_1\le \alpha=2-a-c$ implies $a\le b$, a contradiction. Hence we must have $\alpha=|b-d|$. Then $2-b-d\le \beta_1\le \alpha=|b-d|$ implies that $b=1$ or $d=1$. $b<a$ implies that $b=1$ is impossible. So $d=1$. Now $2-b-c\le \beta_1\le \alpha=1-b$ implies $c\ge 1$ and hence $c=1$, i.e. we have $c=d=1$. Similarly $d<c$ implies that $a=b=1$. This proves our assertion. \\
Next from the definitions we see that $\beta_1\le \gamma_1$ if and only if all conditions (22)-(24) below are satisfied: \begin{equation}d\ge 0 \mbox{ or } c\le a+b+d \mbox{ or } a+b+2d\ge 0\end{equation}
\begin{equation} b\ge 0 \mbox{ or } a\le c+b+d \mbox{ or } c+d+2b\ge 0\end{equation}
\begin{equation} c\le d+2b \mbox{ or } a\le b+2d \mbox{ or } b+d\ge 0\end{equation}
Note that both conditions $a=b=1$ and $c=d=1$ satisfy the above conditions, so we may exclude these cases in the following argument. First we prove that if $\beta_1=\min\{\alpha,\beta_1,\gamma_1\}$, then $b\ge 0$ or $d\ge 0$. Suppose not, i.e. $b<0$ and $d<0$. Since $\beta_1\le \alpha$, we have $c\le d$ and $a\le b$. So $a+b+2d\ge 0\Rightarrow a+b+d\ge -d\ge 0\ge c$ since $d<0$. Similarly $c+d+2b\ge 0 \Rightarrow a\le b+c+d$. Thus conditions (22) and (23) implies that $b+d\ge c-a$ and $b+d\ge a-c$, i.e. $b+d\ge |a-c|\ge 0$, a contradiction. So $b\ge 0$ or $d\ge 0$. Next assume that $a\le b,c\le d, b\ge 0$ and $c>a+b+d$. If $d<0$ then condition (22) implies that $a+b+d\ge -d\ge 0\ge c$, a contradiction. So we must have $d\ge 0$. Now if $a>b+c+d$, then $a>b+(a+b+d)+d=a+2(b+d)$ implies that  $b+d<0$, which contradicts $b+d\ge 0$. Thus $a\le b+c+d$. This proves that either (5) or (6) is true. Since $a\le b$ and $c\le a+b+d$ imply that $c\le 2b+d$, we see that condition (5) implies (22)-(24). Similarly condition (6) implies (22)-(24). This concludes the proof for conditions (5)-(8). \\
The proof for conditions (9)-(12) is similar, replacing $b$ by $-b$ and $d$ by $-d$ throughout in the above argument.\\
It follows from conditions (5)-(6) that if $\beta_1=\min(\alpha,\beta_1,\gamma_1)$, then $2-b-d \le \max(2-b-c, 2-a-d)$, so $\beta_1=\beta$. Similarly if $\gamma_1=\min(\alpha,\beta_1,\gamma_1)$, then $2+b+d \le \max(2+b-c, 2-a+d)$, so $\gamma_1=\gamma$.\\
If $C_0C_1=AC_1,C_1C_2,\cdots,C_{n-1}B=C_{n-1}C_n$ is a path on $X$ joining $A,B$, some of the segments $C_{i-1}C_i$ lying in the half space $y\le -x$, we can reflect those segments back to the set $y\ge -x$ through the isometry $f(x,y,z)=(-y,-x,z)$. The resulting union of segments will have the same total length, which by triangle inequality will be shorter than a path in $y\ge -x$ joining $A, B$. So we only have to consider paths in $X\cap\{(x,y,z)\in X: y\ge -x\}$. Then a simple argument, using triangle inequality, proves that shortest path joining $A,B$ must be one of the three types considered above. \\
Next we show that there is a planar minimal path joining $A,B$. Clearly, we only have to consider cases when $\beta$ or $\gamma$ is minimum. First assume that $\beta$ is minimum and $2-a_z-b_x\ge 2-a_y-b_z$, i.e. $a_z+b_x\le a_y+b_z$. By conditions (5)-(8), we also have $a_y\le a_z$ and $b_x\le b_z$. As in Case (a) above,  minimal paths are of the form $AC_1,C_1C_2,C_2B$ where $C_1=(1,y,1), C_2=(x,1,1)$ and $(x,y)$ satisfies \[1+b_x-b_z\le x\le 1+a_y-a_z\] \[\max(a_y+a_z-1,x)\le y\le 1+a_y-a_z.\]
The path is planar if and only if \[(1-a_z)(1-y)(x-b_x)-(y-a_y)(1-x)(1-b_z)=0.\]
If $a_y=a_z$ and $b_z=1$, then $x=b_x, y=1$ satisfy the above three conditions. Otherwise $a_z-a_y+1-b_z>0$ and direct checking shows that for all $x$ in the interval (which is nonempty since $1+b_x-b_z\le A\le 1+a_y-a_z$)\[1-b_x-b_z\le x\le A\] and \[y=\frac{a_y(1-x)(1-b_z)+(1-a_z)(x-b_x)}{(1-x)(1-b_z)+(1-a_z)(x-b_x)}\] the above three conditions are satisfied, and the corresponding path is planar; here \[A=\frac{b_x(a_z-a_y)+1-b_z}{a_z-a_y+1-b_z}.\]

\begin{co}\rm Let $A, B$ be defined as in the above theorem. There is a shortest path of at most two legs joining $A$ and $B$ if and only if at least one of the conditions (1)-(4) is satisfied.\end{co}
\begin{theo}\rm Let $n=3$. Let $A=(1,a,b)\in X, B=(-1,c,d)\in X$. Then the geodesic distance between $A, B$ is equal to the minimum of the following twelve numbers, $s_1$ to $s_{12}$:
\begin{eqnarray}s_1 &=& 4-a-c\\
s_2 &=& 4+a+c\\
s_3 &=& 4-b-d\\
s_4 &=& 4+b+d\\
s_5 &=& \max\{2-a-d, 4-b-c\}\\
s_6 &=& \max\{2-a+d, 4+b-c\}\\
s_7 &=& \max\{2+a-d, 4-b+c\}\\
s_8 &=& \max\{2+a+d, 4+b+c\}\\
s_9 &=& \max\{2-b-c, 4-a-d\}\\
s_{10} &=& \max\{2+b-c, 4-a+d\}\\
s_{11} &=& \max\{2-b+c, 4+a-d\}\\
s_{12} &=& \max\{2+b+c, 4+a+d\}
\end{eqnarray}
The minimum is $s_1$ if and only if at least one of the following conditions is satisfied:
\begin{eqnarray}|b|\le a &\mbox{ and }& |d|\le c\\
c=-d=1  &\mbox{ and }& b\ge -a\\
c=d=1  &\mbox{ and }& b\le a\\
a=-b=1  &\mbox{ and }& d\ge -c\\
a=b=1  &\mbox{ and }& d\le c
\end{eqnarray}
 The minimum is $s_6$ if and only if at least one of the following conditions is satisfied: (an "and" should substitute for each comma in each condition)
\begin{eqnarray}& & a+b\le 0, c+d\ge 0, |a+d|\le c-b, b+d\le a+c, b\le 0, c\ge 0\\
& & b=-1, c\ge 0, -c\le d\le 1+a+c\\
& & b=-1, d\ge 1-a-c\\
& & c=1 , b\le 0 ,b+d-1\le a\le -b\\
& & c=1 , a\le -1-b-d
\end{eqnarray}
\end{theo}
{\bf Proof.} It is clear that $s_1$ is the shortest length of paths of the form $A, C_1, C_2, B$, where $C_1=(1,1,z_1), C_2=(-1,1,z_2)$ for some $z_1,z_2\in [-1,1]$, and the shortest paths of this type have $|z_1-b|\le 1-a, |z_2-d|\le 1-c, -1\le z_1\le 1, -1\le z_2\le 1$; $s_2$ is the shortest length of paths of the form $A, C_1, C_2, B$, where $C_1=(1,-1,z_1), C_2=(-1,-1,z_2)$ for some $z_1,z_2\in [-1,1]$, and the shortest paths of this type have $|z_1-b|\le 1+a, |z_2-d|\le 1+c, -1\le z_1\le 1, -1\le z_2\le 1$; $s_3$ is the shortest length of paths of the form $A, C_1, C_2, B$, where $C_1=(1,y_1,1), C_2=(-1,y_2,1)$ for some $y_1,y_2\in [-1,1]$, and the shortest paths of this type have $|y_1-a|\le 1-b, |y_2-c|\le 1-d, -1\le z_1\le 1, -1\le z_2\le 1$; $s_4$ is the shortest length of paths of the form $A, C_1, C_2, B$, where $C_1=(1,y_1,-1), C_2=(-1,y_2,-1)$ for some $y_1,y_2\in [-1,1]$, and the shortest paths of this type have $|y_1-a|\le 1+b, |y_2-c|\le 1+d, -1\le z_1\le 1, -1\le z_2\le 1$.\\ \\
Consider paths of the form $A,C_1,C_2,C_3,B$, where $C_1=(1,y,-1),C_2=(x,1,-1), C_3=(-1,1,z)$ for some $x,y,z\in [-1,1]$.
Since $C_1, C_3$ satisfy (9), the shortest distance between $C_1$ and $C_3$ (for given $y,z$) is $\max\{2, 2-y+z\}$ by Theorem 1. It follows that the length of shortest paths of the form  is the minimum of \[g(y,z):=\max(|y-a|,1+b)+\max(2, 2-y+z)+\max(1-c,|z-d|)\]
for $y,z\in [-1,1]$. Since $1+b+2+1-c=4+b-c$ and $y-a+2-y+z+d-z=2-a+d$, the minimum of $g$ is greater than or equal to $s_6$. Now we show that it attains $s_6$. \\
Case (a): $s_6=2-a+d$. This is equivalent to $1+a+b\le c+d-1$. Since $-1\le 1+a+b$ and $c+d-1\le 1$, one sees that $g(y,z)=s_6$ if and only if \[ 1+a+b\le y\le z\le c+d-1;\] and such $(y,z)$ exist. From the proof of Theorem 1, for the path to be shortest, $-1\le x\le \min(y,z)=y.$\\
Case (b): $s_6=4+b-c$. This is equivalent to $c+d-1\le 1+a+b$. It follows that $g(y,z)=s_6$ if and only if $y,z\in [-1,1]$ and $c+d-1\le z\le y\le 1+a+b$. $(y,z)$ satisfies the condition if and only if \[\max(c+d-1,-1)\le z\le y\le \min(1+a+b,1);\] such $(y,z)$ exist for we may take  $y=z=\max(c+d-1,-1)$. From the proof of Theorem 1, for the path to be shortest, $z\le x\le y.$\\
This completes the proof that $s_6$ is the length of the shortest path of the above mentioned form joining $A,B$. \\
The proofs for other $s_i,i\ge 5$ are similar. For example, $s_5$ is the length of shortest path of the form $A,C_1,C_2,C_3,B$, where $C_1=(1,y,-1), C_2=(x,-1,-1), C_3=(-1,-1,z)$ for some $x,y,z\in [-1,1]$; $s_{12}$ is the length of shortest path of the form $A,C_1,C_2,C_3,B$, where $C_1=(1,1,z), C_2=(x,1,1), C_3=(-1,y,1)$ for some $x,y,z\in [-1,1]$.\\
Next we prove that $s_1\le s_j$ for all $j\neq 1$ if and only if $|b|\le a \mbox{ and } |d|\le c$. By comparing $s_1$ with other $s_j$ one finds that $s_1$ is minimal if and only if all the conditions below must be true.
\begin{eqnarray} a+c &\ge& |b+d|\\
 a\ge b &\mbox{or}& c=-d=1\\
  a\ge -b &\mbox{or}& c=d=1\\
   2+d\le 2a+c &\mbox{or}&  a+b+2c\ge 0\\
 2\le 2a+c+d &\mbox{or}& b\ge a+2c\\
  c\ge d &\mbox{or}& a=-b=1\\
 c+d\ge 0 &\mbox{or}& a=b=1\\
 2+b\le 2c+a &\mbox{or}& d\le 2a+c\\
 2\le a+b+2c &\mbox{or}& 2a+c+d\ge 0
 \end{eqnarray}
 Suppose all conditions (47)-(55) are satisfied. If $c=1 \mbox{ and } d=-1$, then (47) implies that $b\ge -a$ and then the remaining conditions (48)-(55) follows. Thus we have condition (37). Conditions (38)-(41) are similarly obtained. If none of the conditions $c=-d=1,c=d=1,a=-b=1$ and $a=b=1$ are true, then conditions (48),(49),(52) and (53) imply that $|b|\le a$ and $|d|\le c$. Then $a\ge 0, c\ge 0$ and we see that (47)-(55) hold. This proves (37)-(41).\\ \\
 By comparing $s_6$ with other $s_j$'s, we get that $s_6$ is minimal if and only if all the conditions below are satisfied:
 \begin{eqnarray}& & a+b\le 0\\
 & & b\le a+2c \mbox{ and } d\le 2+2a+c\\
 & & 2b+d\le c \mbox{ and } b+2d\le 2+a\\
 & & c+d\ge 0\\
 & & (2+a+b+d\le c \mbox{ or }b\le 0)\mbox{ and } (b+c+d\le 2+a \mbox{ or } d\le 0)\\
 & & (b\le c \mbox{ or } 2+b+d\le a+c) \mbox{ and } (b+d\le 2+a+c \mbox{ or } d\le a)\\
 & & (c\ge 0\mbox{ or } 2+b\le a+c+d)\mbox{ and } (d\le 2+a+b+c \mbox{ or } a\ge 0)\\
 & & b=-1 \mbox{ or } a+b+d\le c\\
 & &  a+b\le c+d\\
 & & b+d\le a+c\\
 & & b\le a+c+d\mbox{ or } c=1
 \end{eqnarray}
(65) is the result of comparing $s_6$ with $s_{11}$; note that $1+b\le c \mbox{ and } d\le 1+a$ implies  $b+d\le a+c$. Since (56) and (59) implies (64), (64) can be dropped from the list.\\ (65) implies that $2+2a+c=2+a+a+c\ge 2+a+b+d\ge d$; (65) and (59) implies $a+2c=c+a+c\ge b+d+c\ge b$. Similarly (65) and (56) implies (58). Thus (57) and (58) can be dropped from the list.\\
 Assume $b=-1$ in (63). Then (56),(60),(61) are true. (59) implies (66) since $-b+a+c+d=1+a+c+d\ge 0$.  (65) and (62) implies that $c\ge 0$ or $d\ge 1-a-c$.  So conditions (56)-(66) reduce to (59): $d\ge -c$, (62): $c\ge 0$ or $d\ge 1-a-c$, and (65): $d\le 1+a+c$. If $d\ge 1-a-c$ then $d\ge -c$ (since $1-a\ge 0$) and $1+a+c\ge 2-d\ge d$. So conditions (59) and (65) follows from the condition $d\ge 1-a-c$. This proves (43) and (44). \\
Similarly (45), (46) follow if $c=1$ in condition (66). \\
Excluding $b=-1$ or $c=1$, condition (63) and (66) is equivalent to \begin{equation} |a+d|\le c-b.\end{equation}
Condition (67) implies $b\le c$, and condition (65) obviously implies $b+d\le 2+a+c$. Thus (67) and (65) implies (61), which can then be dropped from the list.\\
Condition (65) implies $d\le a+c-b\le c+b+a+2$, and $b+c+d\le a+2c\le 2+a$. With (67), $2+b\le a+c+d \Rightarrow c\ge 2+b-a-d\ge 2+b+b-c \Rightarrow c\ge 1+b\ge 0$. With (67), $2+a+b+d\le c \Rightarrow b\le c-a-d-2\le c+c-b-2 \Rightarrow b\le c-1\le 0$. This proves that the conditions (56),(59),(65) and (67) change (60) to $b\le 0$ and (62) to $c\ge 0$. This proves (42). \\
We need to prove that other paths do not yield minimal length......

\begin{rem}\rm By symmetry, conditions for other $s_j$'s  being minimal are obtained by appropriately replacing the letters in (37)-(41) or (42)-(46). For example, replacing a by $-b$, b by $-a$, c by $-d$ and d by $-c$ in (37)-(41), we get that $s_4$ is minimal if and only if at least one of the following condition is satisfied:
\begin{eqnarray}|a|\le -b &\mbox{ and }& |c|\le -d\\
c=-d=1  &\mbox{ and }& b\le -a\\
c=d=-1  &\mbox{ and }& b\le a\\
a=-b=1  &\mbox{ and }& d\ge -c\\
a=b=-1  &\mbox{ and }& d\le c
\end{eqnarray}
Replacing b by $-b$, d by $-d$ and keeping a, c unchanged, we get that $s_5$ if and only if at least one of the following conditions is satisfied: (an "and" should substitute for each comma in each condition)
\begin{eqnarray}& & a\le b, d\le c, |a-d|\le c+b, a+b+c+d\ge 0, b\ge 0, c\ge 0\\
& & b=1, c\ge 0, -1-a-c\le d\le c\\
& & b=1, d\le -1+a+c\\
& & c=1 , b\ge 0 ,-b-d-1\le a\le b\\
& & c=1 , a\le -1+b+d
\end{eqnarray}
\end{rem}
\begin{rem}\rm Examples of $A=(1,a,b), B=(-1,c,d)$ where $s_i<s_j$ for $i\neq j$, $i=1,2,\cdots,12$, are \\ $[i,a,b,c,d]=[1,1/20,0,1/20,0],[2,-1,-19/20,-1,-19/20],
[3,-19/20,1,-19/20,1],$ $[4,-19/20,-1,-19/20,-1], [5,-19/20,1/20,19/20,0], [6,-19/20,-1,1/20,0],$\\ $[7,-19/20,1,-1,-19/20], [8, -19/20,-1,-1,-19/20], [9,1/20,-19/20,-1/20,1],$\\ $[10, 1/20,0,-19/20,-1], [11, -1,-19/20,-19/20,1],[12,-1,-19/20,-19/20,-1],$ respectively. In the above example for $s_6$, $2-a+d=4+b-c$. An example of $s_6<s_j, j\neq 6$ and $2-a+d<4+b-c$ is $[-19/20,-19/20,1/20,0]$. An example of $s_6<s_j, j\neq 6$ and $2-a+d>4+b-c$ is $[-19/20,-1,1/20,1/20]$.\\
$A=(1,1/30,-1), B=(-1,-1/30,1)$ is an example which satisfies (44), but fails (42). $A=(1,1/30,-1), B=(-1,0,1)$ is an example which satisfies (43), but fails (42).
\end{rem}
\begin{rem}\rm Note that if $|a|\le -b \mbox{ and } |d|\le c$ then condition (42) is satisfied and $s_6$ is minimal.
\end{rem}
\begin{theo}\rm Let $n\ge 3$. Let $A=(1,a_2,a_3,\cdots,a_n)\in X, B=(b_1,1,b_3,\cdots,b_n)\in X$. For any subset $S$ of $\{i: 3\le i\le n\}$, let $L(S)=\{|a_i-b_i|: 3\le i\le n, i\notin S\}$. Then the geodesic distance between $A, B$ is equal to the minimum of the following numbers:($i,j,k,l,\cdots $ are distinct numbers in $\{3,4,\cdots,n\}$.)
\begin{eqnarray}& &\max\{2-a_2-b_1, |a_3-b_3|,\cdots, |a_n-b_n|\}\\
& &\max\{2-a_2\pm b_i, 2\pm a_i-b_1, L(i)\}\\
& &\max\{2-a_2\pm b_i, 2\pm a_j-b_1, 2\pm a_i\pm b_j, L(i,j)\}\\
& &\max\{2-a_2\pm b_i, 2\pm a_j-b_1, 2\pm a_k\pm b_j, 2\pm a_i\pm b_k, L(i,j,k)\}\\
& & \cdots\\
& &\max\{2-a_2\pm b_i, 2\pm a_j-b_1,2\pm a_k\pm b_j,\cdots, 2\pm a_i\pm b_l\}
\end{eqnarray}
where in (83) the subscript of $a$ in the term preceding $2\pm a_i\pm b_l$ is $l$, and $(i,j,k,\cdots,l)$ is a permutation of $(3,\cdots,n)$.
\end{theo}
{\bf Proof.} An explanation of (79)-(83) is necessary: for a fixed $i$, (79) should have been two numbers, one being
\[\max\{2-a_2+ b_i, 2+ a_i-b_1, L(i)\},\] and the other
\[\max\{2-a_2- b_i, 2- a_i-b_1, L(i)\}.\]
(The choice of signs for $a_i$ and $b_i$ must be the same.) Similarly for fixed $i,j,i\neq j$,(80) should have been four numbers:
\[\max\{2-a_2+ b_i, 2+ a_j-b_1, 2+ a_i+ b_j, L(i,j)\}\]
\[\max\{2-a_2+ b_i, 2- a_j-b_1, 2+ a_i- b_j, L(i,j)\}\]
\[\max\{2-a_2- b_i, 2+ a_j-b_1, 2- a_i+ b_j, L(i,j)\}\]
\[\max\{2-a_2- b_i, 2- a_j-b_1, 2- a_i- b_j, L(i,j)\}.\]
(The choice for the sign of $b_i$ is the same as for $a_i$, and the choice for the sign of $b_j$ is the same as for $a_j$.) Same explanation applies to (81)-(83). The total number of terms in (78)-(83) with various $i,j,k,\cdots$ is $[(n-2)!2^{n-2}\sqrt{e}]$, where $[\cdot]$ stands for the integer part (see Remark below).\\
Each number in (78)-(83) is the length of a shortest path of a certain form. For the case $n= 5, i=3,j=5$, \\$\max\{2-a_2+b_i, 2-a_j-b_1, 2+a_i-b_j, L(i,j)\}$ is the length of the shortest path of the form
$AC_1,C_1C_2, C_2C_3, C_3B$, where $C_1=(1,y,z,u,1), C_2=(x,y_1,-1,u_1,1), C_3=(x_1,1,-1,u_2,v),$  for some $x,x_1,y,y_1,z,u,u_1,u_2,v$ in $[-1,1]$. \\
Generally, $\max\{2-a_2\pm b_i, 2\pm a_j-b_1,2\pm a_k\pm b_j,\cdots, 2\pm a_i\pm b_l,L(i,j,k,\cdots,l)\}$ is the length of a shortest path of the form $AC_1,C_1,C_2,\cdots,C_mB$, given schematically in the following table. Here $m$ is one more than the number of integers in $\{i,j,k,\cdots,l\}$, and $i,j,k,\cdots,l$ need not be in increasing order.
\[\begin{array}{cccccccccccccc}A:& 1&a_2&a_3&\cdots&a_i&\cdots&a_j&\cdots&a_k&\cdots&a_l&\cdots&a_n\\
C_1:&1&&&&&&\mp 1&&&&&&\\
C_2:&&&&&&&\mp 1&&\mp 1&&&&\\
&&&&\mp 1&&&&&\mp 1&&&&\\
\cdot&\cdot&\cdot&\cdot&\cdot&\cdot&\cdot&\cdot&\cdot&\cdot&\cdot&\cdot&\cdot&\cdot\\
\cdot&\cdot&\cdot&\cdot&\cdot&\cdot&\cdot&\cdot&\cdot&\cdot&\cdot&\cdot&\cdot&\cdot\\
\cdot&\cdot&\cdot&\cdot&\cdot&\cdot&\cdot&\cdot&\cdot&\cdot&\cdot&\cdot&\cdot&\cdot\\
\cdot&&&&&\mp 1&&&&&&\mp 1&\\
C_m:&&1&&&\mp 1&&&&&&&&\\
B:&b_1&1&b_3&\cdots&b_i&\cdots&b_j&\cdots&b_k&\cdots&b_l&\cdots&b_n
\end{array}\]

\begin{rem}\rm Let $m$ be a positive integer. Then \begin{eqnarray*}& & 1+2m+2^2m(m-1)+\cdots +2^mm!\\
&=&m!2^m\left(1+\cdots +\frac{2^{i-m}}{(m-i)!}+\cdots +\frac{2^{-m}}{m!}\right)\\
&=&m!2^m\left(\sqrt{e}-R_m\right)
\end{eqnarray*} where $R_m=\sum_{i=m+1}^\infty (.5)^i/i!$. Note that $m!2^mR_m<1$, so \\$1+2m+2^2m(m-1)+\cdots +2^mm!=[m!2^{m}\sqrt{e}]$. \end{rem}
\begin{theo}\rm Let $n\ge 3$. Let $A=(1,a_2,a_3,\cdots,a_n)\in X, B=(-1,b_2,b_3,\cdots,b_n)\in X$.  Then the geodesic distance between $A, B$ is equal to the minimum of the following numbers:($i,j,k,l,m,\cdots $ are distinct numbers $\{2,3,\cdots,n\}$.)
\begin{eqnarray}& &4\pm a_i\pm b_i\\
& &\max\{2\pm a_i\pm b_j, 4\pm a_j\pm b_i\}\\
& &\max\{2\pm a_i\pm b_j, 2\pm a_k\pm b_i, 4\pm a_j\pm b_k\}\\
& & \cdots\\
& &\max\{2\pm a_i\pm b_j, 2\pm a_k\pm b_i, 2\pm a_l\pm b_k,\cdots, 4\pm a_j\pm b_m\}
\end{eqnarray}
where in (88) the subscript of $a$ in the term preceding $4\pm a_j\pm b_m$ is $m$, and if (88) is the last one on the list then $(i,j,k,l,\cdots,m)$ is a permutation of $(2,\cdots,n)$.
\end{theo}
{\bf Proof.} As explained in the proof of the preceding, (84)-(88) consist of $[(n-1)!2^{n-1}\sqrt{e}]-1$ terms.\\
For $n=6, i=3,j=5,$ $\max\{a+a_i-b_j, 4-a_j+b_i\}$ is the length of a shortest path of the form $AC_1,C_1C_2,C_2C_3,C_3B$, where $C_1=(1,y,z,u,1,w), C_2=(x,y_1,-1,u_1,1,w_1), C_3=(-1,y_2,-1,u_2,v,w_2)$ for some $x,y,z,u,v,w,y_1,y_2,u_1,u_2,w_1,w_2\in [-1,1]$.\\
Generally, $\max\{2\pm a_i\pm b_j, 2\pm a_k\pm b_i,2\pm a_l\pm b_k,\cdots, 4\pm a_j\pm b_m\}$ is the length of a shortest path of the form $AC_1,C_1,C_2,\cdots,C_pB$, given schematically in the following table. Here $p$ is one more than the number of integers in $\{i,j,k,l\cdots,m\}$, and $i,j,k,l\cdots,m$ need not be in increasing order.
\[\begin{array}{ccccccccccccccc}A:& 1&a_2&\cdots&a_i&\cdots&a_j&\cdots&a_k&\cdots&a_l&\cdots&a_m&\cdots&a_n\\
C_1:&1&&&&&\mp 1&&&&&&&&\\
C_2:&&&&&&\mp 1&&&&&&\mp 1&&\\
\cdot&\cdot&\cdot&\cdot&\cdot&\cdot&\cdot&\cdot&\cdot&\cdot&\cdot&\cdot&\cdot&\cdot&\cdot\\
\cdot&\cdot&\cdot&\cdot&\cdot&\cdot&\cdot&\cdot&\cdot&\cdot&\cdot&\cdot&\cdot&\cdot&\cdot\\
\cdot&\cdot&\cdot&\cdot&\cdot&\cdot&\cdot&\cdot&\cdot&\cdot&\cdot&\cdot&\cdot&\cdot&\cdot\\
C_{p-2}:&&&&&&&&\mp 1&&\mp 1&&&&\\
C_{p-1}:&&&&\mp 1&&&&\mp 1&&&&&&\\
C_p:&-1&&&\mp 1&&&&&&&&&&\\
B:& -1&b_2&\cdots&b_i&\cdots&b_j&\cdots&b_k&\cdots&b_l&\cdots&b_m&\cdots&b_n
\end{array}\]

\end{document}